\documentclass[12pt]{article}

\usepackage{hyperref,yfonts,pifont,mathrsfs,latexsym,palatino,amsthm,amssymb,amsmath,amsfonts,graphicx,color,tikz}

\allowdisplaybreaks

\makeatletter
\newcommand{\leqnomode}{\tagsleft@true\let\veqno\@@leqno}
\newcommand{\reqnomode}{\tagsleft@false\let\veqno\@@eqno}
\makeatother

\begin{document}

\section*{Cantor's Powerset Theorem, Graph--Theoretically}
\vspace{1cm}

\begin{center}
\textbf{Saeed Salehi}\footnote{%
Harish and Bina Shah School of AI and Computer Science, Plaksha University, IT City Road, Sector 101A, Mohali, Punjab 140306, India. \textbf{\texttt{root@SaeedSalehi.ir}}}
\end{center}

\noindent 
A directed graph is a pair $G\!=\!(V,E)$, where $V$ is a nonempty set of vertices and $E$ is a binary relation of edges between the vertices, $E\!\subseteq\!V^2$. Here, $vEw$ denotes the fact that the vertex $v\!\in\!V$ is  connected by an edge in $E$ to the vertex $w\!\in\!V$, and 
$\boldsymbol\neg vEw$
indicates that $v$ is not
$E$-connected to $w$. The outgoing set of a vertex $v\!\in\!V$ is the set of all the vertices to which $v$ is $E$-connected, i.e., ${\sf Out}(v)\!=\!\{w\!\in\!V\mid vEw\}$.

If a combinatorist, or an AI tool, is asked {\em why there should exist a set of vertices
that is not the outgoing set of any vertex}, the most likely answer will be that {\em the subsets of $V$ outnumber all the outgoing sets}, which are at most as many as the number of all vertices. Indeed, if $V$ has $n$ elements, then there are $2^n$ subsets of $V$, and there are at most $n$ outgoing sets, and since $2^n\!>\!n$, then, by the pigeonhole principle, there should exist a set of vertices that is not the outgoing set of any vertex. If the graph is infinite, then the argument remains the same; just the ``cardinality'' replaces the ``number of elements.''

The ingenious diagonal argument of {\sc Cantor}  (1845--1918)
explicitly constructs,  with an elementary proof, a set of vertices that is not equal to any outgoing set.
A vertex $v\!\in\!V$ is called  looped  when it is connected to itself by an edge in $E$, i.e., when  $vEv$; otherwise it is  called  unlooped (when 
$\boldsymbol\neg vEv$). {\sc Cantor}'s (anti-)diagonal set is the set of all unlooped vertices,
$$\mathcal{D}\!=\!\{u\!\in\!V\mid 
\boldsymbol\neg uEu
\},$$ which he showed to be unequal to any outgoing set. Here is a proof.

\begin{quote}
\begin{itemize}
\item[]
  \textsf{\small
  Fix a vertex $v\!\in\!V$. If $vEv$, then $v\!\in\!{\sf Out}(v)$ and $v\!\not\in\!\mathcal{D}$, so we have $v\!\in\!{\sf Out}(v)\!\setminus\!\mathcal{D}$.  If
$\boldsymbol\neg vEv$,
then $v\!\not\in\!{\sf Out}(v)$ and  $v\!\in\!\mathcal{D}$, so we have $v\!\in\!\mathcal{D}\!\setminus\!{\sf Out}(v)$.
Thus, we always have $v\!\in\![{\sf Out}(v)\!\setminus\!\mathcal{D}]\cup[\mathcal{D}\!\setminus\!{\sf Out}(v)]$. Therefore, we have $\mathcal{D}\!\neq\!{\sf Out}(v)$   for each and every vertex $v\!\in\!V$.
    }
\end{itemize}
\end{quote}

There are some other variants of this  marvelous proof. For a positive integer $n$, let $\mathcal{D}_n$ be the set of all the vertices that are not connected to themselves by an $E$-path of length $n\!+\!1$, i.e.,
$$\mathcal{D}_n\!=\!\{v\!\in\!V\mid \textrm{there are } \textbf{no }
w_1,\cdots,w_n \textrm{ with } vEw_1E\cdots Ew_nEv\}.$$
Let $\mathcal{D}_\infty$ be the set of all the vertices for which there is no infinite path starting from them, i.e.,
$$\mathcal{D}_\infty\!=\!\{v\!\in\!V\mid \textrm{there are } \textbf{no }  w_1,\cdots,w_n,\cdots \textrm{ with } vEw_1E\cdots Ew_nE\cdots\}.$$
The proofs are similar to {\sc Cantor}'s \textbf{\cite{Cantor}}; the idea of $\mathcal{D}_n$ is from {\sc Quine}  \textbf{\cite[{\rm \S24}]{Quine}}, and $\mathcal{D}_\infty$ is {\sc Raja}'s \textbf{\cite{Raja}}.
The following argument seems to be new (see \textbf{\cite{Salehi1}} and \textbf{\cite{Salehi2}} for some other proofs of {\sc Cantor}'s powerset theorem). Fix $S$ to be a nonempty  set of natural numbers, i.e., $\emptyset\!\neq\!S\!\subseteq\!\{0,1,2,3,\cdots\}$.  Let $\mathcal{D}_S$ be the set of all vertices that are not connected to themselves by an $E$-path of length $n\!+\!1$ for some $n\!\in\!S$, i.e.,
$$\mathcal{D}_S\!=\!\{v\!\in\!V\mid \textrm{there are } \textbf{no }
w_1,\cdots,w_n \textrm{ with } 0\!<\!n\!\in\!S \textrm{ and } $$
$$\qquad \qquad \qquad \qquad \qquad \qquad
vEw_1E\cdots Ew_nEv; \textrm{ also }
\boldsymbol\neg vEv
\textrm{ when } 0\!\in\!S\}.$$

\bigskip

\noindent
We now prove that none of
$\mathcal{D}_n/\mathcal{D}_\infty/\mathcal{D}_S$ can be the outgoing set of a vertex.

\begin{quote}
\begin{itemize}
\item[]
  \textsf{\small
For a fixed $v\!\in\!V$, we show that $\mathcal{D}_x\!\subseteq\!{\sf Out}(v)$ implies ${\sf Out}(v)\!\not\subseteq\!\mathcal{D}_x$ ($x\!\in\!\{n,\infty,S\}$).
If $vEv$ holds, then $v\!\in\!{\sf Out}(v)$, but $v\!\not\in\!\mathcal{D}_x$  since  $vEvE\cdots EvE\cdots$. Thus, $v\!\in\!{\sf Out}(v)\!\setminus\!\mathcal{D}_x$, so ${\sf Out}(v)\!\not\subseteq\!\mathcal{D}_x$.
Let us now assume that 
$\boldsymbol\neg vEv$
holds. Thus,   $v\!\not\in\!{\sf Out}(v)$, so $v\!\not\in\!\mathcal{D}_x$ by the assumption $\mathcal{D}_x\!\subseteq\!{\sf Out}(v)$.
Then, for some $w_1,\cdots,w_k,\cdots$ we have   $vEw_1E\cdots Ew_kE\cdots$. Now,  $w_1\!\in\!{\sf Out}(v)$, but $w_1\!\not\in\!\mathcal{D}_x$ since
$w_1E\cdots Ew_nE\cdots$ in the case of $x\!=\!\infty$, and
$w_1E\cdots Ew_nEvEw_1$
in the case of $x\!\in\!\{n,S\}$. Thus, $w_1\!\in\!{\sf Out}(v)\!\setminus\!\mathcal{D}_x$, so ${\sf Out}(v)\!\not\subseteq\!\mathcal{D}_x$.}
\end{itemize}
\end{quote}

\noindent
Let us notice that
$\mathcal{D}_\infty\!\subseteq\!\mathcal{D}_n\!\subseteq\!\mathcal{D}$ for every $n\!>\!0$, and that $\mathcal{D}_S\!=\!\bigcap_{n\in S}\mathcal{D}_n$
 with the convention $\mathcal{D}_0\!=\!\mathcal{D}$.


\begin{thebibliography}{99}


\bibitem{Cantor}
{\sc G. Cantor}, {\em \"Uber eine elementare Frage der Mannigfaltigkeitslehre},
\textbf{\textit{Jahresbericht der Deutschen Mathematiker-Vereinigung}} 1 (1891), 75--78.
(English translation: ``On an elementary question in the theory of manifolds'', W. Ewald, ed. (1996), {\sl From  Kant to  Hilbert, Vol.\ 2}, {\sc isbn}:~\href{https://isbnsearch.org/isbn/9780198505365}{\tt 9780198505365}, Oxford University Press, pp. 920--922.)



\bibitem{Quine}
{\sc W. Quine}, \textbf{\em Mathematical Logic}, Harvard University Press (revised 1951).
{\sc isbn}:~\href{https://isbnsearch.org/isbn/9780674554504}{\tt 9780674554504}



\bibitem{Raja}
{\sc N. Raja}, {\em A negation-free proof of Cantor's theorem}, \textbf{\textit{Notre Dame Journal of  Formal  Logic}} 46:2 (2005), 231--233.
{\sc doi}: \href{https://doi.org/10.1305/ndjfl/1117755152}{\tt 10.1305/ndjfl/1117755152}



\bibitem{Salehi1}
{\sc S. Salehi}, {\em A non-constructive proof of Cantor’s theorem},
arXiv:\href{https://doi.org/10.48550/arXiv.2510.14534}{\tt 2510.14534} [math.LO] (2025), pp. 1--2.



\bibitem{Salehi2}
{\sc S. Salehi}, {\em Cantor's non-equinumerosity theorems, inductively},
arXiv:\href{https://doi.org/10.48550/arXiv.2510.15321}{\tt 2510.15321} [math.LO] (2025), pp. 1--7.


\end{thebibliography}
\end{document}